\documentclass[12pt]{article}
\usepackage{amsmath,amssymb}
\usepackage{tikz}
\usepackage{color}
\usepackage{caption}

\newtheorem{theorem}{Theorem}[section]

\newcommand{\proof}{\noindent{\bf Proof.\ }}
\newcommand{\qed}{\hfill $\square$ \bigskip}

\usepackage{graphicx}
\begin{document}

\title{Characterization of classes of graphs with large  general  position number}

\author{
	Elias John Thomas $^{a}$
	\and
	Ullas Chandran S. V. $^{b}$
}

\maketitle
\begin{center}
	$^a$ Department of Mathematics, Mar Ivanios College, Thiruvananthapuram-695015, Kerala, India;\quad 
	{\tt eliasjohnkalarickal@gmail.com}
	\medskip
	
	$^b$ Department of Mathematics, Mahatma Gandhi College, Kesavadasapuram,  Thiruvananthapuram-695004, Kerala, India;\quad 
	{\tt svuc.math@gmail.com} 
	\medskip		
\end{center}

\maketitle

\begin{abstract}
 Getting inspired by the famous no-three-in-line problem and by the general position subset selection problem from discrete geometry, the same is introduced into graph theory as follows. A set $S$ of vertices in a graph $G$ is a general position set if no element of $S$ lies on a geodesic between any two other elements of $S$. The cardinality of a largest general position set is the general position number ${\rm gp}(G)$ of $G.$  In \cite{ullas-2016} graphs $G$ of order $n$ with ${\rm gp}(G)$ $\in \{2, n, n-1\}$ were characterized.
In this paper, we characterize the classes of all connected graphs of order $n\geq 4$ with the general position number $n-2.$
\end{abstract}

\medskip\noindent

{\bf AMS Subj.\ Class:} 05C12; 05C69.
\section{Introduction}
\label{sec:intro}

 The general position problem in graphs was introduced by P.~Manuel and S.~Klav\v{z}ar \cite{manuel-2018} as a natural extension of the well known century old Dudeney's no-three-in-line problem and the general position subset selection problem from discrete geometry~\cite{dudeney-1917,froese-2017,payne-2013}. The general position problem in graph theory was introduced in~\cite{manuel-2018} as follows. A set $S$ of vertices in a graph $G$ is a \emph{ general position set}  if no element of $S$ lies on a geodesic between any two other elements of $S$. A largest general position set is called a $gp$-$set$ and its size is the\emph{ general position number}(gp-number, in short), ${\rm gp}(G)$, of $G.$

The same concept was in use two years earlier in~\cite{ullas-2016} under the name {\em geodetic irredundant  sets}. The concept was defined in a different method, see the preliminaries below. In~\cite{ullas-2016} it is proved that for a connected graph of order $n,$ the complete graph of order $n$ is the only graph with the largest general position number $n$; and ${\rm gp}(G)= n-1$ if and only if $G=K_1+\bigcup_j m_jK_{j}$ with $\sum m_j\geq 2$ or $G=K_n-\{e_1, e_2, \dots, e_k\}$ with $1\leq k\leq n-2$, where $e_i$'s all are edges in $K_n$ which are incident to a common vertex $v$. In~\cite{manuel-2018}, certain general upper and lower bounds on the gp-number are proved. In the same paper it is proved that  the general position problem is NP-complete for arbitarary graphs.  The $\rm gp$-number for a large class of subgraphs of the infinite grid graph, for the infinite diagonal grid, and for Bene\v{s} networks were obtained in the subsequent paper~\cite{P manuel-2018}. Anand et al.\cite{elias-2019} gives a characterization of general position sets in arbitrary graphs. As a consequence, the $\rm gp$-number of graphs of diameter $2$, cographs, graphs with at least one universal vertex, bipartite graphs and their complements were obtained. Subsequently, $\rm gp$-number for the complements of trees, of grids, and of hypercubes were deduced  in  \cite{elias-2019}. Recently, in \cite {ghorbani-2019} a sharp lower bound on the gp-number is proved for Cartesian products of graphs. In the same paper the $\rm gp$-number for joins of graphs, coronas over graphs, and line graphs of complete graphs are determined. Recent developments on general position number can be seen in \cite {balázs-2019}.
\section{Preliminaries}
\label{sec:preliminary}

 Graphs used in this paper are finite, simple and undirected.  The distance $d_G(u,v)$ between $u$ and $v$ is the minimum length of an $u,v$-path. An $u,v$-path of  minimum length is also called an $u,v$-{\it geodesic}.
  The maximum distance between all pairs of vertices of $G$ is the {\em diameter},  ${\rm diam}(G)$, of $G$. A subgraph $H$ of a graph $G$ is \emph{isometric subgraph} if $d_H(u, v) = d_G(u, v)$ for all $u, v \in
  V(H)$. A The {\em interval} $I_G[u,v]$ between vertices $u$ and $v$ of a graph $G$ is the set of vertices that lie on some $u,v$-geodesic of $G$. For $S\subseteq V(G)$ we set $I_G[S]=\bigcup_{_{u,v\in S}}I_G[u,v]$. We may simplify the above notation by omitting the index $G$ whenever $G$ is clear from the context.

A set of vertices $S\subseteq V(G)$ is a {\em general position set} of $G$ if no three vertices of $S$ lie on a common geodesic in $G$. A gp-set is thus a largest general position set. Call a vertex $v\in T\subseteq V(G)$ to be an {\em interior vertex} of $T$, if $v \in I [T - \{v\}]$. Now, $T$ is a general position set if and only if $T$ contains no interior vertices. In this way general position sets were introduced in~\cite{ullas-2016} under the name geodetic irredundant sets. The maximum order of a complete subgraph of a graph $G$ is denoted by $\omega(G).$ Let $\eta(G)$ be the maximum order of an induced complete multipartite subgraph of the complement of $G.$ Finally, for $n\in \mathbb{N}$ we will use the notation $[n] = \{1,\ldots,n\}$.

In this paper, we make use of the following results.
\begin{theorem}
{\rm\cite{ullas-2016}}
\label{thm1}
Let $G$ be a connected graph of order $n$ and diameter $d$. Then ${\rm gp}
(G)\leq n-d+1$.
\end{theorem}

\begin{theorem}
{\rm \cite{ullas-2016}}
\label{thm2}
For any cycle $C_n$ $(n\geq 5)$, ${\rm gp}(C_n)=3$.
\end{theorem}

We recall the characterization of general position sets from \cite{elias-2019}, for which we need some additional information.
Let $G$ be a connected graph, $S\subseteq V(G)$, and ${\cal P} = \{S_1, \ldots, S_p\}$ a partition of $S$. Then ${\cal P}$ is \emph{distance-constant} if for any $i,j\in [p]$, $i\ne j$, the distance $d(u,v)$, where $u\in S_i$ and $v\in S_j$ is independent of the selection of $u$ and $v$.  If ${\cal P}$ is a distance-constant partition, and $i,j\in [p]$, $i\ne j$, then let $d(S_i, S_j)$ be the distance between a vertex from $S_i$ and a vertex from $S_j$. Finally, we say that a distance-constant partition ${\cal P}$ is {\em in-transitive} if $d(S_i, S_k) \ne d(S_i, S_j) + d(S_j,S_k)$ holds for arbitrary pairwise different $i,j,k\in [p]$.
\begin{theorem}
{\rm\cite{elias-2019}}
\label{thm:gpsets}
Let $G$ be a connected graph. Then $S\subseteq V(G)$ is a general position set if and only if the components of $G[S]$ are complete subgraphs, the vertices of which form an in-transitive, distance-constant partition of $S$.
\end{theorem}

\begin{theorem}
{\rm\cite{elias-2019}}
\label{thm:diameter2}
If ${\rm diam}(G) = 2$, then ${\rm gp}(G) = \max\{\omega(G), \eta(G)\}$.
\end{theorem}

\section{The characterization}
\label{sec:characterization}
In the following, we characterize all connected graphs $G$ of  order $n\geq 4$ with the $\rm gp$- number $n-2.$ Since the complete graph $K_n$ is the only connected graph of  order $n$ with the $\rm gp$-number $ n$, by Theorem \ref {thm1}, we need to consider only graphs with diameter 2 or 3. First, we introduce four families of graphs with the diameter 3; and four families of graphs with the diameter 2.

Let ${\cal F}_1$ be the collection of  all graphs obtained from the cycle $C:u_1, u_2, u_3, u_4,u_1$  by adding $k$ new vertices $v_1, v_2, \dots, v_k$($k\geq 1)$ and  joining each $v_i, i \in [k]$ to the vertex $u_1.$  Graphs from the family ${\cal F}_1$ are presented in Figure 1.
\begin{figure} 
\begin{center}
    \includegraphics[scale=0.1]{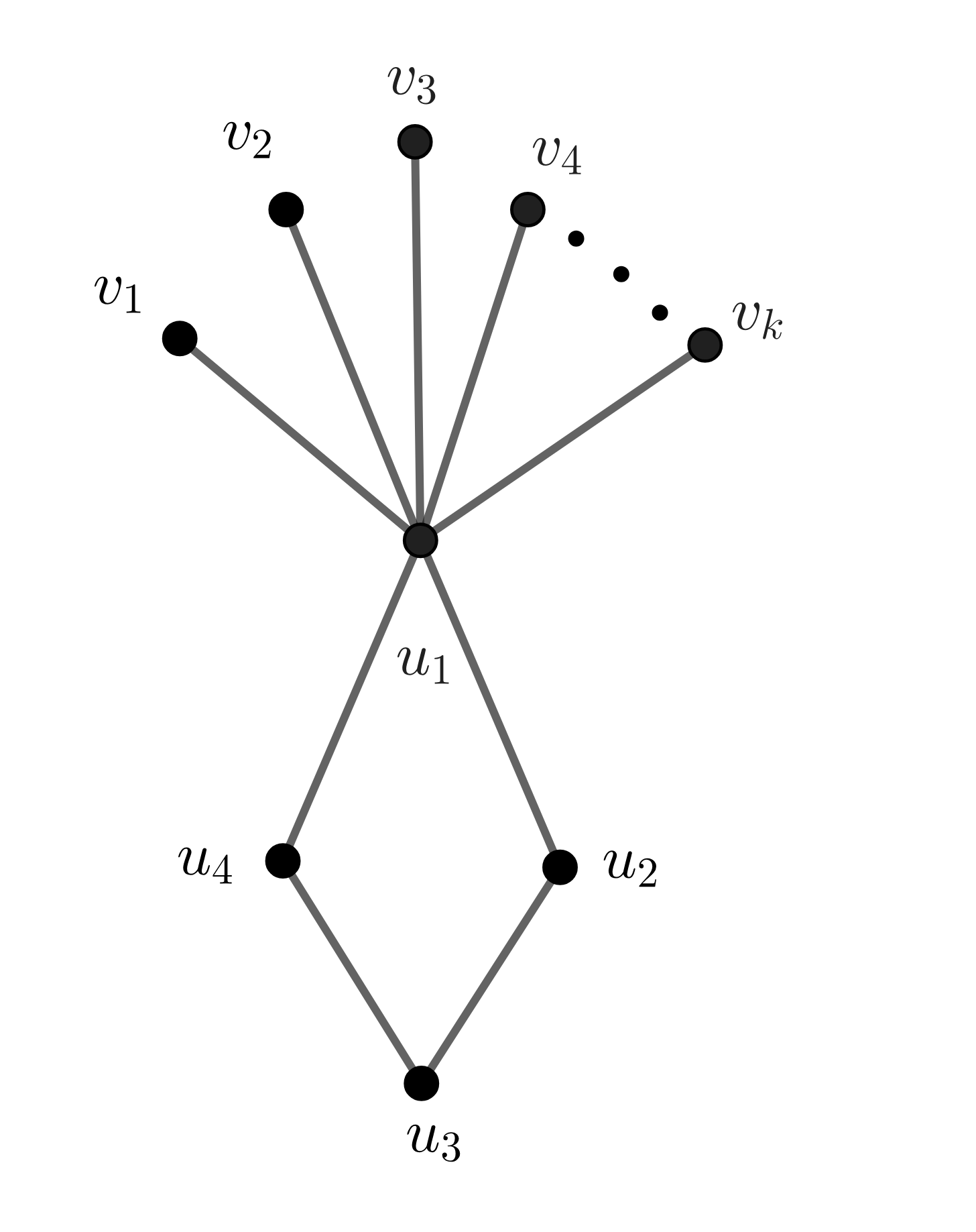}
    \captionof{figure}{Family ${\cal F}_1$}
\end{center}
\end{figure}

  Let ${\cal F}_2$ be the collection of  all  graphs obtained from the path $P_2:x,y$  and complete graphs $K_{n_1 }, K_{n_2},\dots, K_{n_r }(r\geq 1),$  $ K_{m_1 }, K_{m_2},\dots,K_{m_s }(s\geq 1)$ and $K_{l_1 }, K_{l_2}, \dots, K_{l_t }$ (possibly complete graphs of this kind may be empty), by joining both $x$ and $y$ to all vertices of $K_{l_1 }, K_{l_2}, \dots, K_{l_t }; $ joining $x$ to all vertices of $K_{n_1 }, K_{n_2},\dots, K_{n_r };$ and joining $y$ to all vertices of $K_{m_1 }, K_{m_2},\dots,K_{m_s }.$ Graphs from the family ${\cal F}_2$ are presented in Figure 2. Trees with diameter 3 are called double stars and they belong to the class ${\cal F}_2.$

\begin{figure} 
\begin{center}
    \includegraphics[scale=0.14]{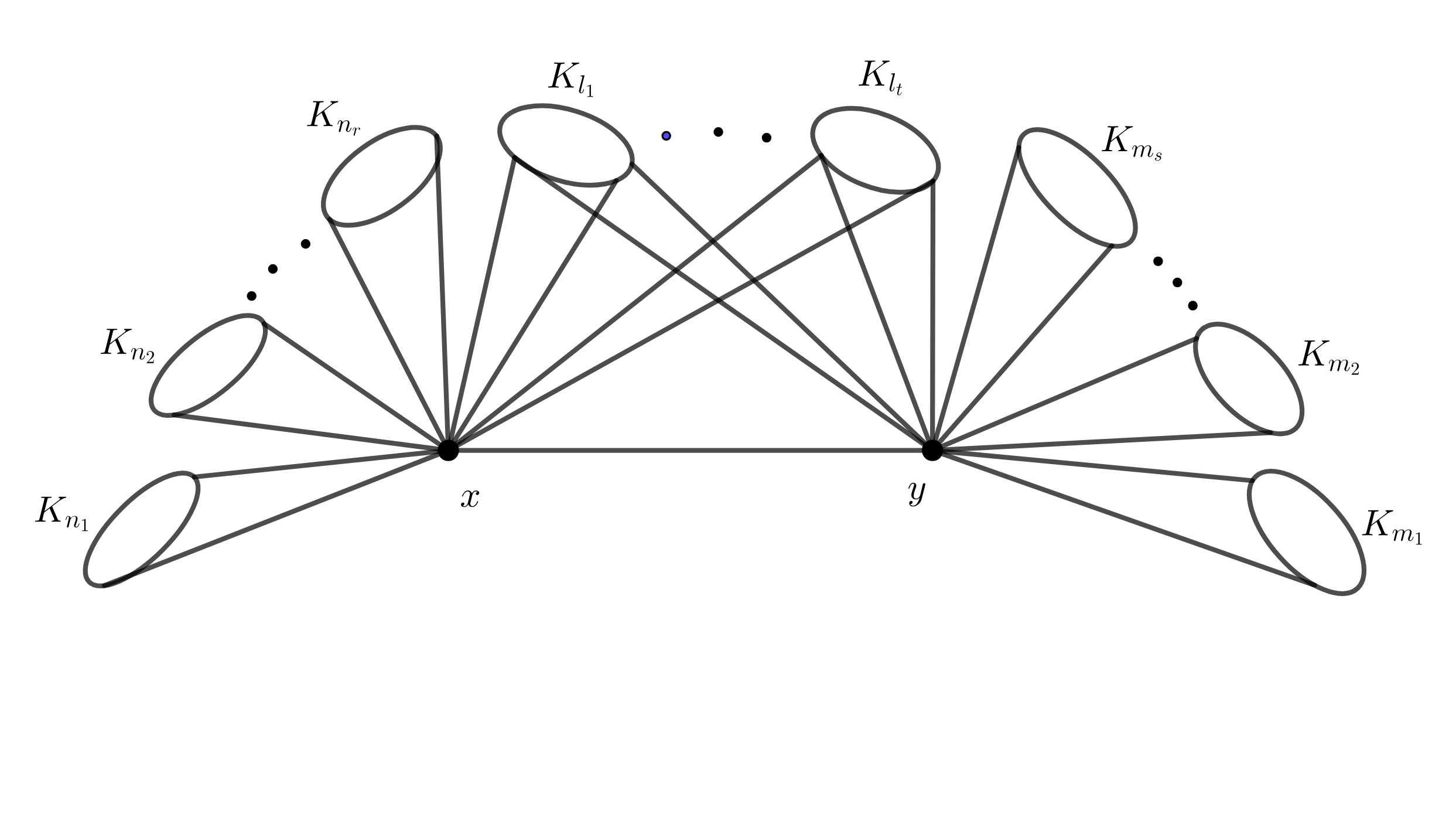}
    \captionof{figure}{Family ${\cal F}_2$}
\end{center}
\end{figure}

Let ${\cal F}_3$ be the collection of  all  graphs obtained from the  path $P_4:u,x,y,v$ and a complete graph $K_r(r\geq 1)$ by joining both $u$ and $x$ to all vertices of $K_r$ and joining $y$ to a subset $S$ of vertices  of $V(K_r)$ (possibly $S$ may be empty or $S=V(K_r)$). Graphs from the family ${\cal F}_3$ are presented in Figure 3.

\begin{figure} 
\begin{center}
    \includegraphics[scale=0.10]{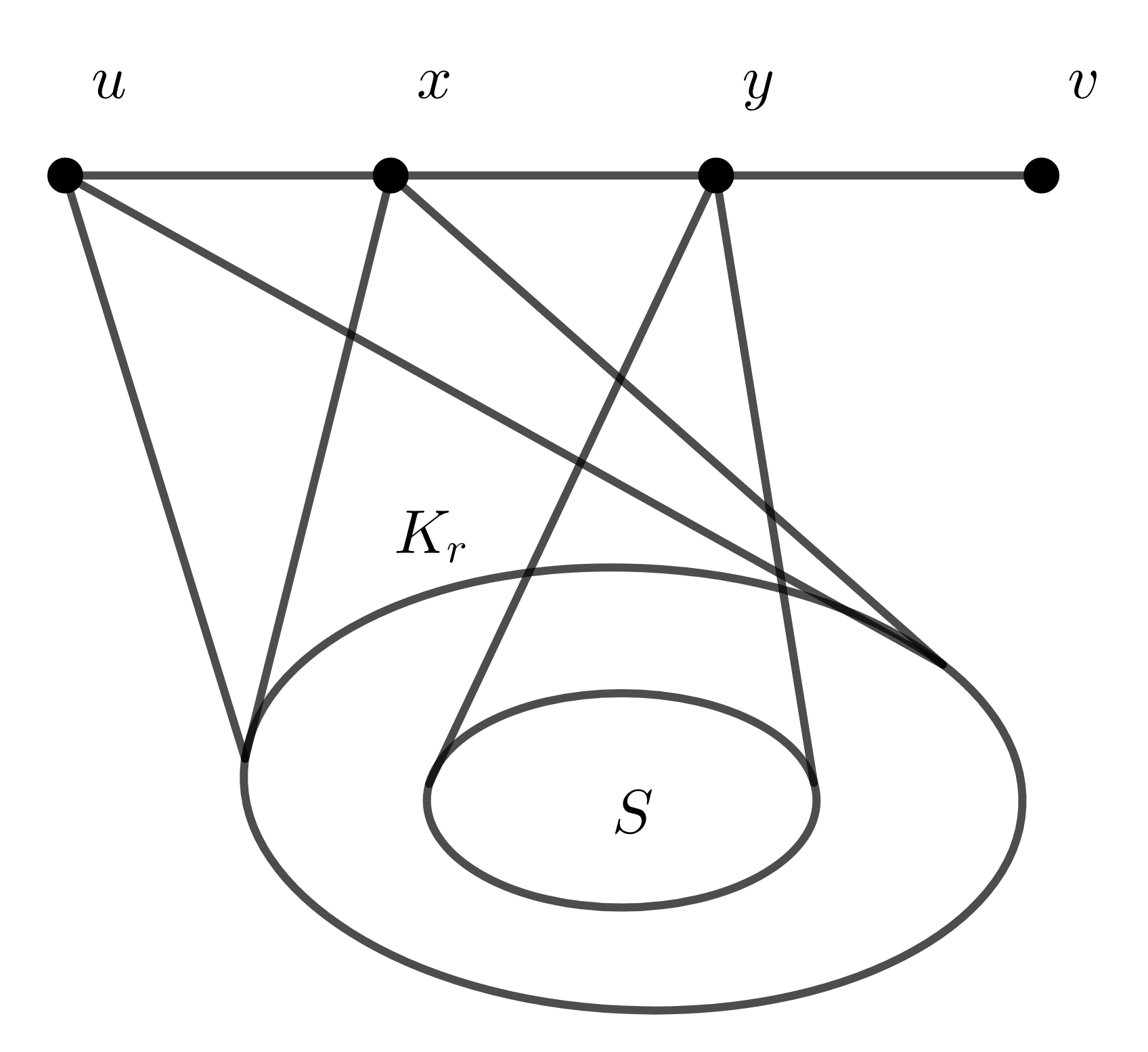}
    \captionof{figure}{Family ${\cal F}_3$}
\end{center}
\end{figure}

Let ${\cal F}_4$ be the collection of  all  graphs obtained from the path $P_3:x,y,v$  and complete graphs $K_q, K_{n_1 }, K_{n_2},\dots, K_{n_r }(r\geq 1),$ $ K_{m_1}, K_{m_2 },\dots,K_{m_s}(s\geq 1)$ by joining $x$ to all vertices of $K_{n_1 }, K_{n_2},\dots, K_{n_r };$  joining $x$ and $v$ to all vertices of $ K_{m_1}, K_{m_2 },\dots,K_{m_s};$ joining $x$ and $y$ to all vertices of $ K_q.$  Graphs from the family ${\cal F}_4$ are presented in Figure 4.

\begin{figure} 
\begin{center}
    \includegraphics[scale=0.15]{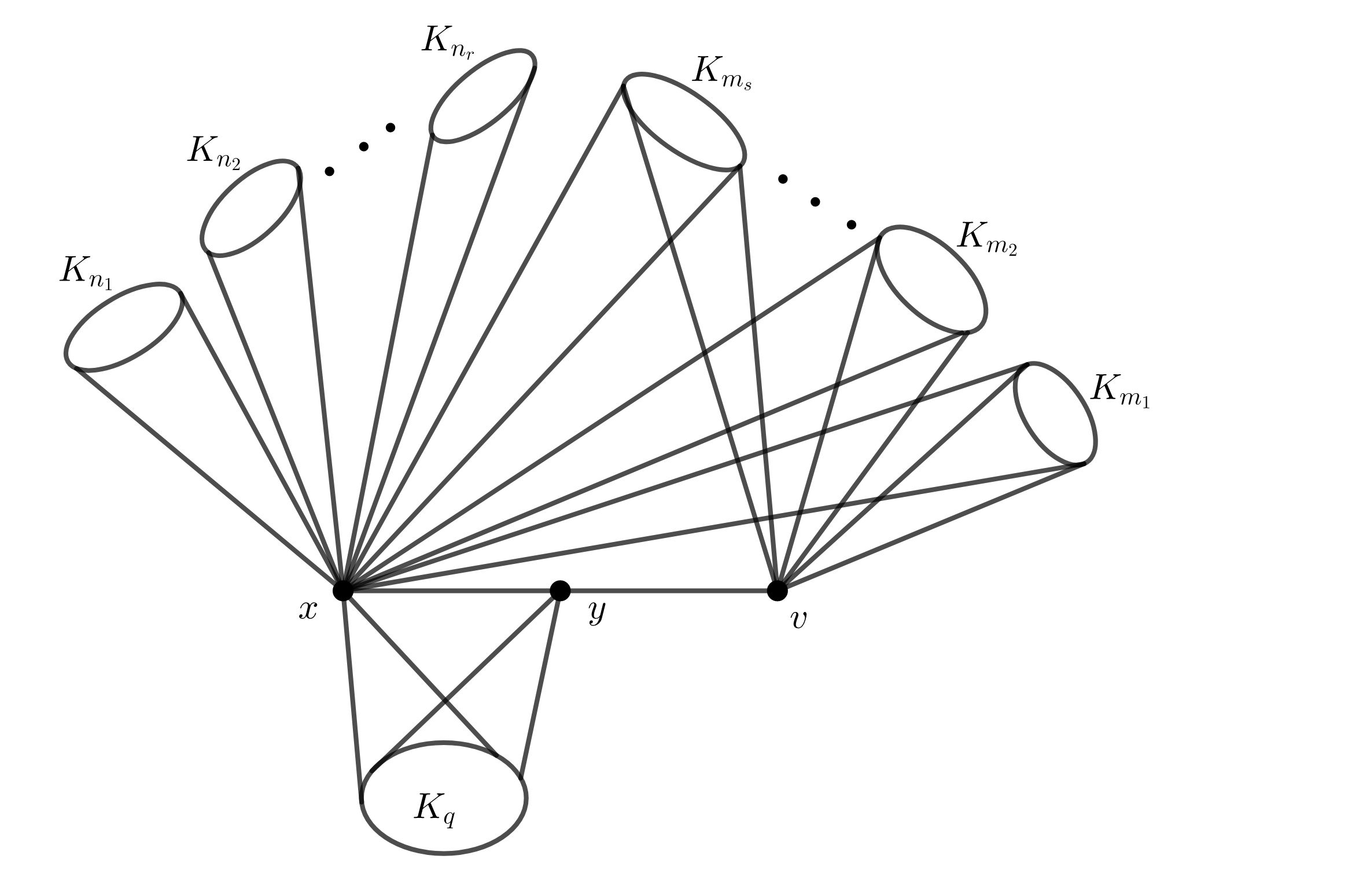}
    \captionof{figure}{Family ${\cal F}_4$}
\end{center}
\end{figure}


 Next, we introduce four families of graphs with diameter 2.

 Let ${\cal F}_5$ be the collection of  all  graphs obtained from the complete graph $K_{n-2} (n\geq 5)$ by adding two new vertices $u$ and $v,$ joining $u$ to all vertices of non-empty subset $S$ of $V(K_{n-2})$ of size at most $n-3$; and joining $v$ to all vertices of non-empty subset $T$ of $V(K_{n-2})$ of size at most $n-3.$ The set $S$ must intersect with the set $T$ so that, the diameter of each graph from the family ${\cal F}_5$ is 2. Graphs from the family ${\cal F}_5$ are presented in Figure 5.

\begin{figure} 
\begin{center}
    \includegraphics[scale=0.11]{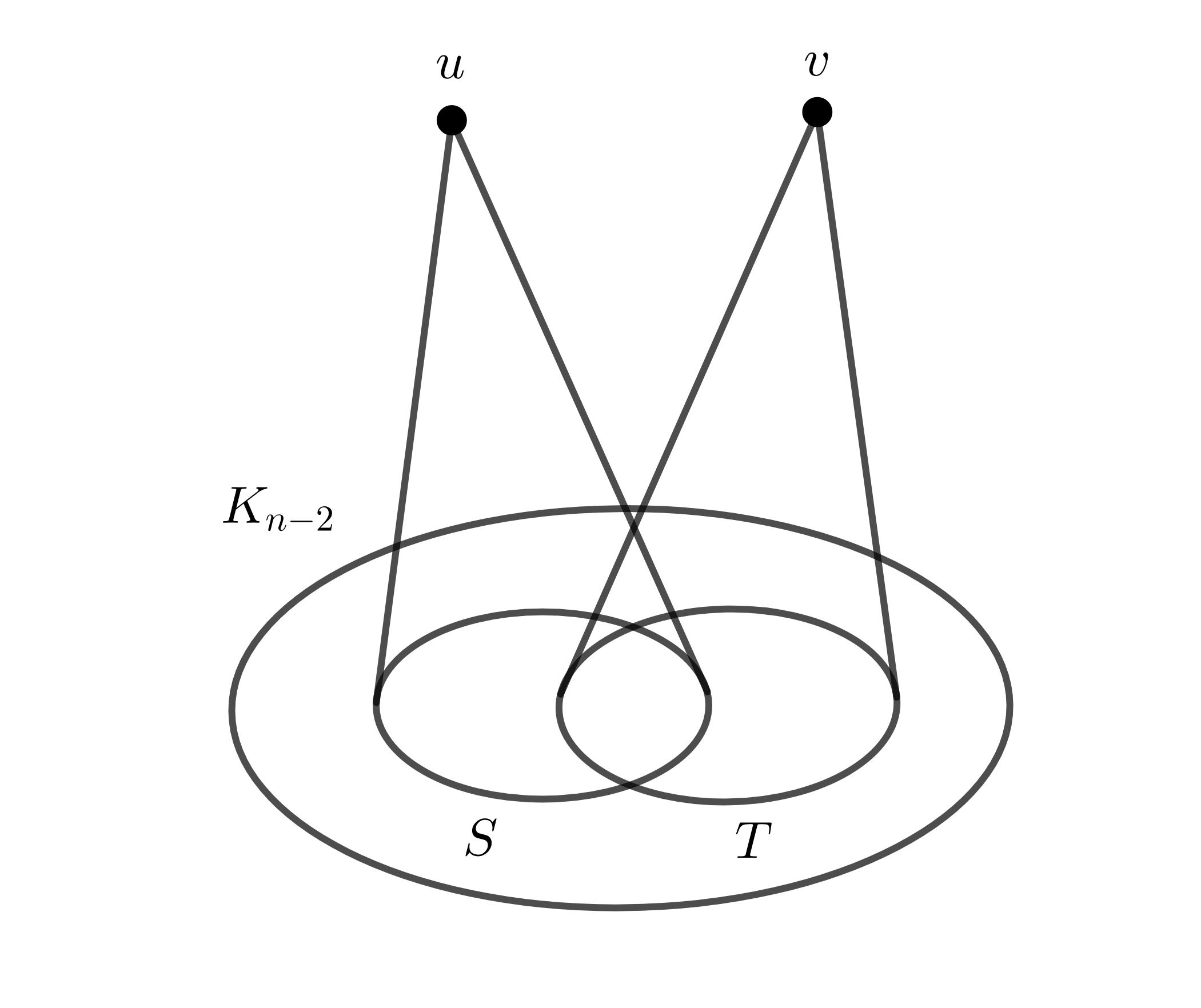}
    \captionof{figure}{Family ${\cal F}_5$}
\end{center}
\end{figure}

Let ${\cal F}_6$ be the collection of  all  graphs obtained from the family ${\cal F}_5$ by adding the edge $uv.$ Moreover; in this case, the set $S$ may be disjoint with the set $T.$ Graphs from the family ${\cal F}_6$ are presented in Figure 6.
\begin{figure} 
\begin{center}
    \includegraphics[scale=0.11]{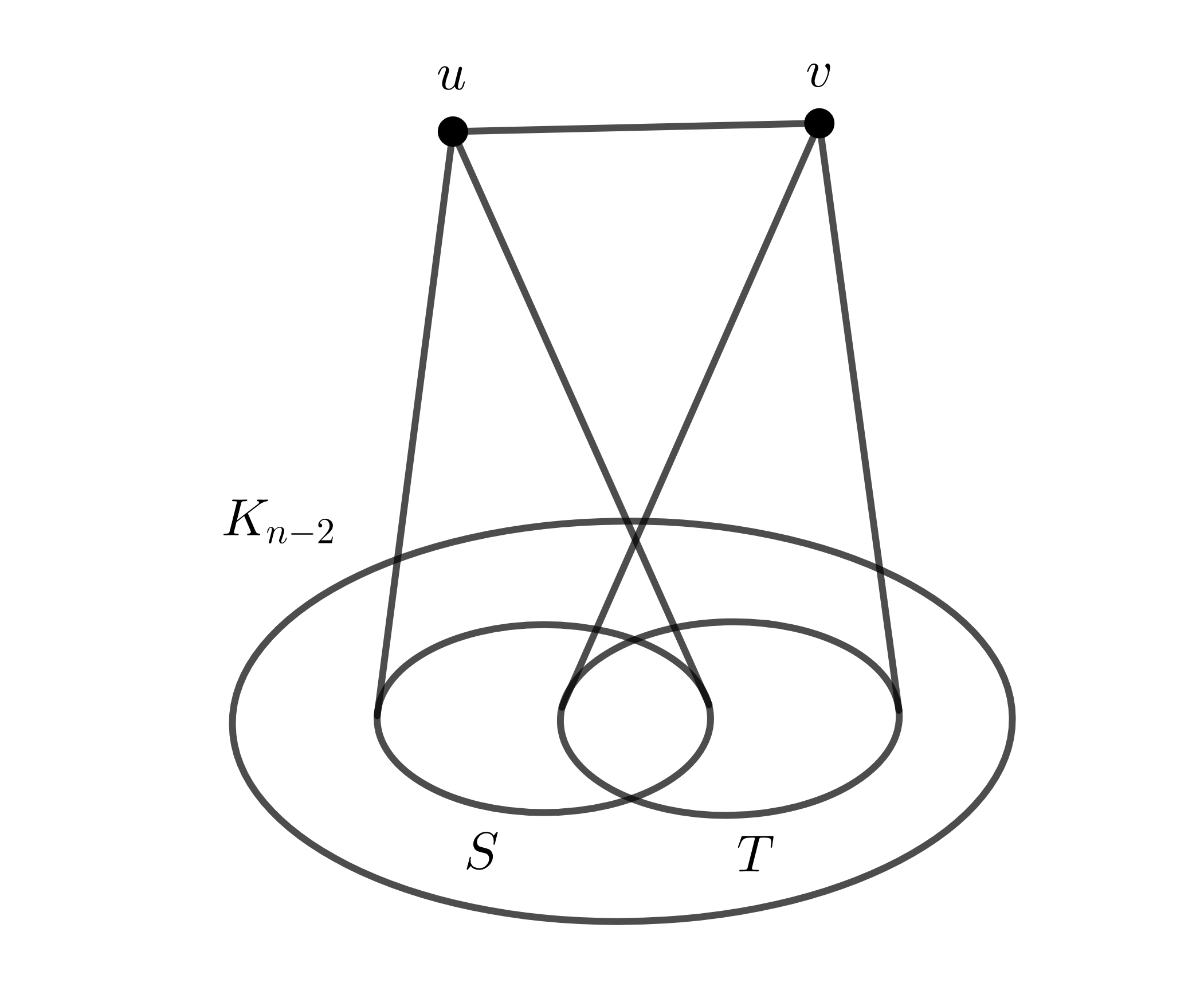}
    \captionof{figure}{Family ${\cal F}_6$}
\end{center}
\end{figure}

Let ${\cal F}_7$ be the collection of  all  graphs obtained from  the complete graphs $K_{n_1 }, K_{n_2},\\\dots, K_{n_r }(r\geq 2)$ by adding two new vertices $x$ and $y,$ joining $x$ to a non-empty subset $S_i$ of $V(K_{n_i})$ for all $i \in [r];$ and $y$ to a non-empty subset $T_i$ of $V(K_{n_i})$ for all $i \in [r]$ (the edges are in  a way that for any $ u\in V(K_{n_i})$ and $v\in V(K_{n_j})$ with $i\neq j$ must have a common neighbor). Moreover, for some $i \in [r]$; the set $S_i$ must intersect with the set $T_i$ so that, the diameter of each graph from the family ${\cal F}_7$ is 2.  Graphs from the family ${\cal F}_7$ are presented in Figure 7. It is clear that both  $C_4$ and $C_5$ belong to class ${\cal F}_7.$

\begin{figure} 
\begin{center}
    \includegraphics[scale=0.15]{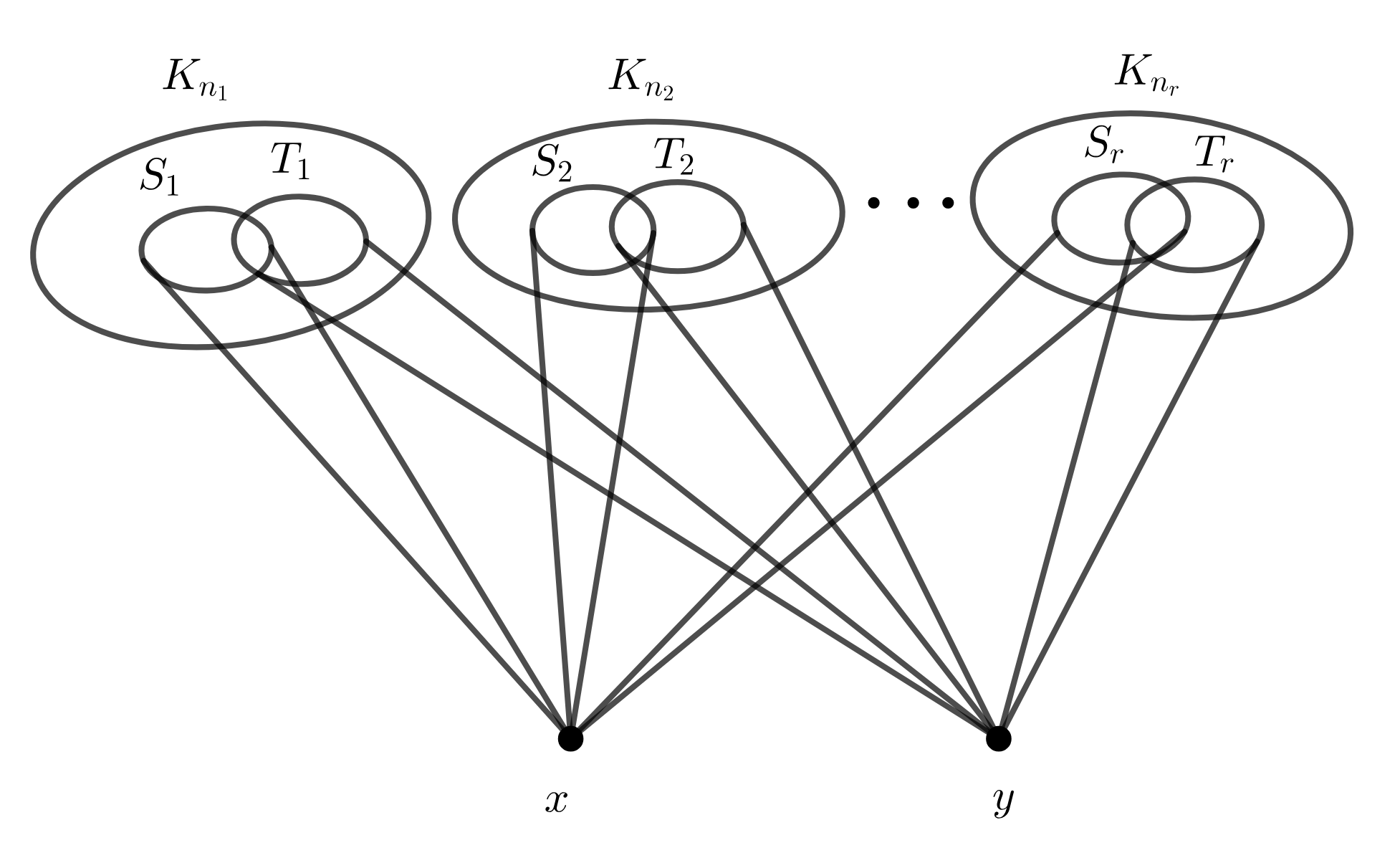}
    \captionof{figure}{Family ${\cal F}_7$}
\end{center}
\end{figure}

Let ${\cal F}_8$ the collection of all graphs obtained from the family ${\cal F}_7$ by adding the edge $xy.$  In this case, the set $S_i$ may be disjoint with the set $T_i$ for all $i \in [r].$ Graphs from the family ${\cal F}_8$ are presented in Figure 8.

\begin{figure} 
\begin{center}
    \includegraphics[scale=0.11]{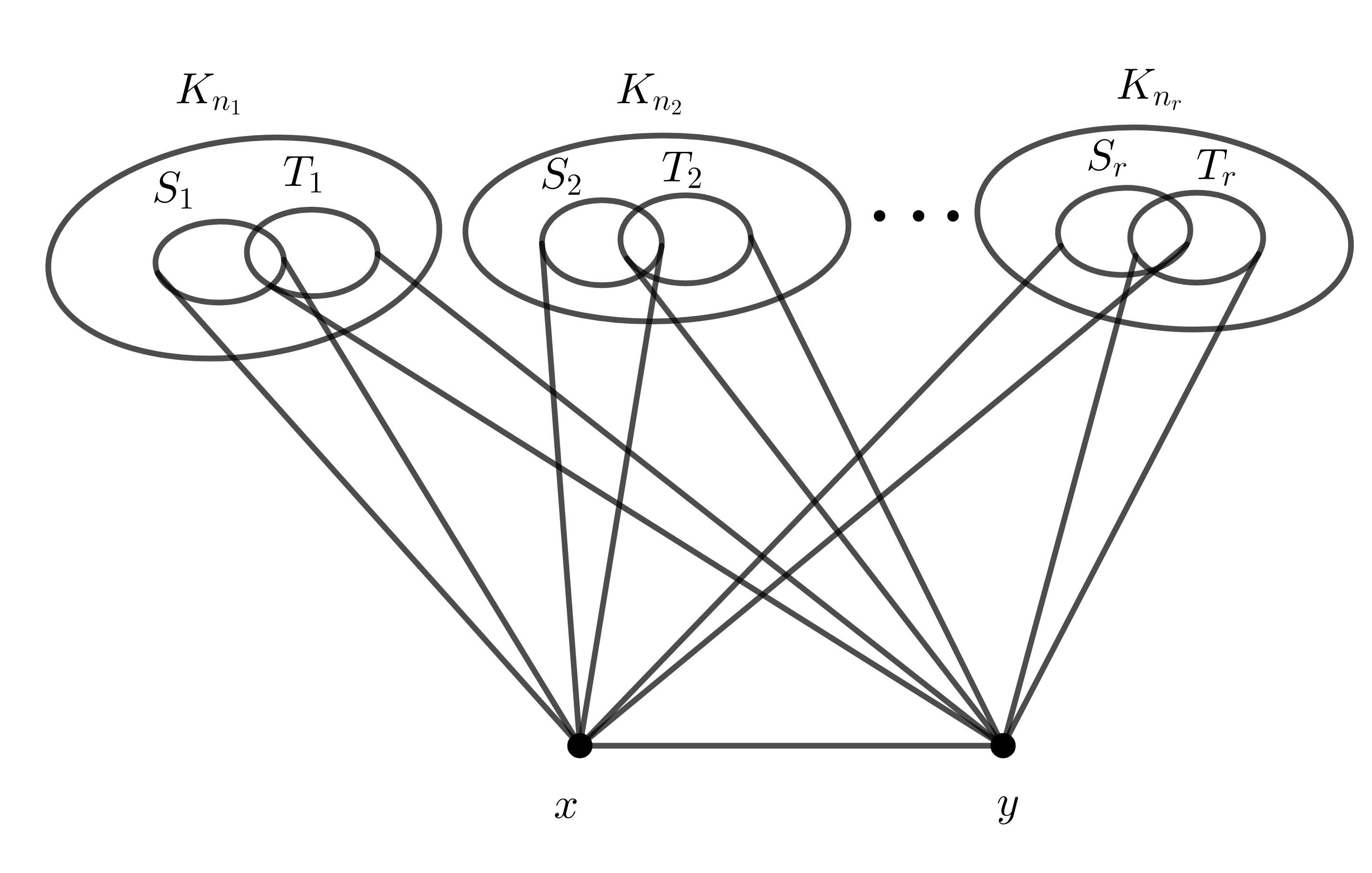}
    \captionof{figure}{Family ${\cal F}_8$}
\end{center}
\end{figure}
\newpage
\begin{theorem}
Let $G$ be a connected graph of order $n\geq 4,$ then ${\rm gp}(G) = n-2$  if and only if $G$ belongs to the family $\cup_{i=1}^8{\cal F}_i.$
\end{theorem}

\rm \proof
First, suppose that $G$ is a connected graph of order $n$ with ${\rm gp}(G)=n-2.$ Then it follows from Theorem \ref{thm1} that ${\rm diam}(G)$ is either 2 or 3. We consider the following two cases.\\
 {\bf Case 1:}
    {\rm diam}($G$) = 3. If $G$ is a tree, then $G$ is a double star and hence it belongs to ${\cal F}_2$.
     So, assume that $G$ has cycles. Let ${\rm girth}(G)$ denotes the length of a  shortest cycle in $G.$ 
     
     Let $C$ be any shortest cycle in $G$. Then it is clear that $C$ is an isometric subgraph of $G$. This shows that if $S$ is a general position set in $G$, then $S\cap V(C)$ is a general position set in $C$. Hence it follows from Theorem \ref {thm2} that any general position set of $G$ contains at most three vertices from the cycle $C.$ Now, since  ${\rm gp}(G)= n-2,$ we have that the length of $C$ is at most 5 and so ${\rm girth}(G)\leq 5.$ 
     
     Next, we claim that there is no connected graph of order $n$ with ${\rm girth}(G) = 5$ and ${\rm gp}(G)= n-2.$ For, assume the contrary that there is a connected graph of order $n$ with ${\rm girth}(G) = 5$ and ${\rm gp}(G)= n-2.$
 Let $C:u_1,u_2,u_3,u_4,u_5,u_1$ be a shortest cycle of length 5 in $G.$  Since ${\rm girth}(G) = 5,$ it follows that the vertices from $N(u_i)$ are independent for all $i \in [5].$ Also, as above we have that any general position set of $G$ has at most three vertices from the cycle $C.$ Let $S$ be a general position set in $G.$ Since ${\rm gp}(G)= n-2,$ we have that $S = V(G)\setminus\{u_i,u_j\}.$ If $u_i$ and $u_j$ are sucessive vertices in $C,$ then it follows that the induced subgraph of $S$ has a $P_3,$ which is impossible. Hence without  loss of generality, we may assume that $i =1$ and $j=3.$ So $S = V(G)\setminus\{u_1,u_3\}.$ Now, since $u_2, u_4, u_5 \in S$ and $N(u_i)$ is independent, by Theorem \ref{thm:gpsets}, it follows that ${\rm deg}(u_i)\leq 3$ for $i = 2, 4, 5.$ Now we claim that ${\rm deg}(u_2)= {\rm deg}(u_4)= {\rm deg}(u_5)= 2.$ Otherwise, we may assume that ${\rm deg}(u_2)=3$ and let $x$ be the neighbour of $u_2$ different from $u_1$ and $u_3.$ Since ${\rm girth}(G) = 5,$ it follows that $x$ is not adjacent with the remaining vertices of $C.$ Now, since $u_2,u_5,x \in S,$ by Theorem  \ref{thm:gpsets}, $d(u_5,x) = d(u_5, u_2)=2.$ Let $P:u_5,y,x$ be a $u_5,x$-geodesic of length 2. Then it is clear that $y\notin V(C)$ and so $y\in S.$ This leads to the fact that induced subgraph of $S$ has a $P_3$, impossible in a general position set. Hence ${\rm deg}(u_2)= 2.$ Similarly $ {\rm deg}(u_4)= {\rm deg}(u_5)= 2.$
 
  Now, if $N(u_1)\neq \emptyset,$ then $u_5 \in I[x, u_4]$ for all $x\in N(u_1)$ (otherwise $S$ contains an induced $P_3$), impossible. Hence $N(u_1)= \emptyset.$ Similarly, $N(u_2)= \emptyset.$ Hence $G\cong C_5.$ But $ {\rm gp}(C_5)= 3=n-2$ and ${\rm diam}(G)= {\rm diam}(C_5)=2.$ Hence there is no connected graph of order $n$ with ${\rm diam}(G) =3,$  ${\rm girth}(G)=5$ and ${\rm gp}(G)= n-2.$ Hence ${\rm girth}(G)$ is at most 4. 
  
  Now, assume that ${\rm girth}(G)=4$ and let  $C:u_1,u_2,u_3,u_4,u_1$ be a shortest cycle of length 4 in $G.$ Since ${\rm diam}(G) = 3,$ we have that $G\ncong C_4.$ Now, we may assume that $u_1 \in V(C)$ be a vertex such that ${\rm deg}(u_1)\geq 3$ and  let $x$ be a neighbour of $u_1$ such that $x\notin V(C).$ Since $S$ is a general position set and $|S| = n-2,$ we have that $S$ contains exactly 2 vertices from $C.$ We claim that $u_1 \notin S.$
  For otherwise assume that $u_1 \in S.$  Since $|S| = n-2$ and $x,u_1 \in S,$ it follows from Theorem \ref{thm:gpsets} that $u_2,u_4 \notin S$ and $u_3 \in S.$ This shows that the path $x,u_1,u_2,u_3$ must be a $x,u_3$- geodesic (otherwise, since $|S|=n-2,$ $S$ contains an induced $P_3.$ Hence $d(x,u_3) \neq d(u_1, u_3),$ which is impossible in a general position set. Hence $u_1 \notin S.$
  
   Now, we claim that $u_1$ is the unique vertex in $C$ with degree at least 3.
   Assume the contrary that there exists $u_j \in C$ with $j\neq 1$ and ${\rm deg}(u_j)\geq 3.$ Then as above we have that $u_j \notin S.$ Now, if $u_i$ and $u_j$ are adjacent vertices in $C,$ then we can assume that $j=2.$
It follows from the fact that $S$ is a general position set of size $n-2,$ $d(u_3,x)=3$ and $u_3,u_4,u_1,x$ is a geodesic in $G,$ where $x$ is a neighbour of $u_1$ such that $x\notin V(C).$ This shows that the vertices $x,u_4,u_3,x$ lie on  a common geodesic, a contradiction. Similarly if $u_1$ and $u_j$ are non adjacent vertices in $C$ then $ u_j = u_3$ and $u_2,u_4$ belong to $S.$ Moreover, as above $S$ is a general position set of size $n-2,$ we have that $x,y\in S$
and $d(x,y) =4,$ where $x\in N(u_1) \setminus V(C)$ and $ y \in N(u_3) \setminus V(C)$, which is impossible. Thus $u_1$ is the unique vertex in $C$ with ${\rm deg}(u_1)\geq 3.$ Also, since ${\rm girth}(G) =4,$ we have that $N(u_1)$ induces an independent set. Hence the graph belongs to ${\cal F}_1.$ 

Now, consider ${\rm girth}(G) =3$ and ${\rm diam}(G)=3.$ Let $P:u,x,y,v$ be a $u,v$- shortest path in $G$ of length 3. Then $S$ contains atmost 2 vertices from $V(P).$ Since $|S| =n-2,$ we have that $S$ contains exactly two vertices from $V(P).$  We consider the following four cases.\\
 {\bf Subcase 1.1:}
 $u,v\in S.$ Then $x,y\notin S.$ Moreover, $S = V(G)\setminus\{x,y\}.$ Now, let $z$ be any neighbour of $u.$ Since $S$ is a general position set of size $n-2,$ it follows that $I[z,v]\subseteq V(P).$  This shows that $d(z,v)\leq 3.$ If $d(z,v)=2,$ then $z$ must be adjacent with $y$ and so $u,z,y,v$ is a $u-v$ geodesic, which contradicts the fact that $S$ is a general position set. Hence $d(z,v)=3$ and since $I[z,v]\subseteq V(P),$ we have that $z$ is adjacent with $x$ but it is not adjacent with $y.$ Similarly, we have that any neighbour of $v$ is adjacent with $y$ but non-adjacent with $x.$ Now, assume that $z$ be any vertex in $G$ such that $z\notin V(P)$ and $z$ is non-adjacent with both $u$ and $v.$ Then as in the previous case, we have that $I[z,v]\subseteq V(P).$ Also, we have $d(z,v) \in \{2,3\}$ and $d(z,u)\in \{2,3\}.$ Hence it follows that $z$ is adjacent to $x$ or $y$ or both. Also, by Theorem \ref{thm:gpsets}, we have that the components of $S$ are in-transitive distance-constant cliques. Hence the graph reduces to the class ${\cal F}_2.$\\
 {\bf Subcase 1.2:}
$u,x\in S.$ Then $y,v\notin S$  and $S = V(G)\setminus \{y,v\}.$ Now, let $z$ be any vertex in $G$ such that $z\notin V(P).$ Then, we have that $I[z,u]\subseteq V(P).$ Moreover, by Theorem \ref{thm:gpsets}, $d(z,u) = d(z,x).$ If $d(z,x)=2,$ then $I[z,x]\subseteq V(P),$ we have that $z$ is adjacent to $y.$ But in this case $d(z,u)$ cannot be equal to 2. Similarly, if $d(z,x)=3$ then $z$ is adjacent with $v$ but not $y.$ Then it is clear that $d(z,u)\neq 3.$ Hence it follows that  $d(z,u) = d(z,x)=1.$ Again by Theorem \ref{thm:gpsets}, $V(G)\setminus\{y,v\}$ induces a clique. Hence the graph reduces to the class ${\cal F}_3.$\\
 {\bf Subcase 1.3:}
  $u,y\in S.$ Then $x,v\notin S$ and $S = V(G)\setminus\{x,v\}.$ Now, for any $z\notin V(P),$ we have that  $I[z,y]\subseteq V(P)$ and  $I[z,u]\subseteq V(P).$ Thus $d(z,y)\leq 3$ for all $z\notin V(P).$ If $d(z,y)=3,$ then $z$ must be adjacent to $u$ and so by Theorem \ref{thm:gpsets},  $d(u,y)=3,$ a contradiction. Thus $d(z,y)\in \{1,2\}.$ If $d(z,y)=1,$ then again by Theorem \ref{thm:gpsets}, we have that $d(u,z)=2$ and so $z$ must be adjacent to $x.$ Moreover, $ \{z\notin V(P): d(z,y)=1\}$ induces a clique. Now, if $d(z,y)=2,$ then by using the same argument, we have that $z$ is either adjacent to $x$ or $z$ is adjacent to both $x$ and $v.$ Hence the graph reduces to class  ${\cal F}_4.$\\
 {\bf Subcase 1.4:}
 $x,y\in S.$ Then $u,v\notin S$ and $S = V(G)\setminus\{u,v\}.$ Now, for any $z\notin V(P),$ as in the previous case we have that  $I[z,x]\subseteq V(P)$ and  $I[z,y]\subseteq V(P).$ Moveover, by Theorem \ref{thm:gpsets},  $d(z,x)=d(z,y).$ Now, if $d(z,x)\neq 1, $ then $d(z,y)\neq 1.$ This shows that $z$ must be  adjacent to both $u$ and $v,$ which is impossible. Hence $d(z,x)=d(z,y) = 1.$ Hence it follows from Theorem \ref{thm:gpsets}, $V(G)\setminus\{u,v\}$ induces a clique. Moreover, since both $x$ and $y$ belong to $S,$ it is clear that $d(u,z)=d(v,z) = 2$ for all $z \notin V(P).$ Hence in this case the graph reduces to the family ${\cal F}_2.$\\
{\bf Case 2:}
${\rm diam}(G)= 2.$ Then by Theorem \ref{thm:diameter2}, we have $ ${\rm gp}($G$)$ = \max\{\omega(G), \eta(G)\}= n-2.$  We consider the following two subcases.\\
\bf{Subcase 2.1:}
$\omega(G)\geq \eta(G).$ \rm Then $ ${\rm gp}($G$)$ =\omega(G)= n-2.$ Let $K$ be a clique of order $n-2$ and let $u,v \in V(G)$ be such that $u,v \notin V(K).$  Then it is clear that $1 \leq {\rm deg}(u) \leq n-3 $ and $1 \leq {\rm deg}(v) \leq n-3. $  Now, if $u$ and $v$ are adjacent in $G,$ then $G$ belongs to the family  ${\cal F}_6.$ Otherwise, $G$ belongs to the family ${\cal F}_5.$\\
 \bf {Subcase 2.2:}
$ \eta(G)>\omega(G).$
\rm Then $ ${\rm gp}($G$)$ =\eta(G)= n-2.$ This shows that the complement of $G$ has  complete mulipartite subgraph $H$ of order $n-2.$ Thus the components of the induced subgraphs of $H$ in $G$ are cliques, say $K_{n_1 }, K_{n_2},\dots, K_{n_r }.$  Moreover $d(u,v)=2$ for all $u \in V(K_{n_i})$ and $v \in V(K_{n_j}).$ Now, let $x$ and $y$ be the vertices in $G$ such that $x, y \notin V(H).$ Then it is clear that  the graph reduces to the family ${\cal F}_8,$ when $x$ and $y$ are adjacent in $G.$ Otherwise it belongs to the family ${\cal F}_7.$

  On the other hand, if $G$ belongs to the family $\cup_{i=1}^8{\cal F}_i$, by Theorems 2.1 and \ref{thm:gpsets}, one can easily verify that ${\rm gp}(G) = n-2.$ This completes the proof.
  
 \qed

\section*{Acknowledgements}

The authors are grateful to the anonymous referees for their valuable suggestions and comments. E.J. acknowledges the University of Kerala for providing JRF for the research work.


\end{document}